\documentclass[12pt]{amsart}
\usepackage{graphicx}
\usepackage{color}
\usepackage{amssymb}
\usepackage{amsmath}
\usepackage{epstopdf}
\usepackage{epsfig}
\usepackage{psfrag}
\numberwithin{figure}{section}
\newtheorem{theorem}{Theorem}[section]
\newtheorem{corollary}[theorem]{Corollary}
\newtheorem{lemma}[theorem]{Lemma}

\theoremstyle{definition}
\newtheorem{definition}[theorem]{Definition}
\newtheorem{remark}[theorem]{Remark}
\title{Stabilizing Heegaard splittings of toroidal 3-manifolds}
\author{Ryan Derby-Talbot}
\address{Department of Mathematics \\ The American University in Cairo}
\email{rdt@aucegypt.edu}

\begin{document}

\maketitle

\newcommand{\hs}{$V \cup_S W$}
\newcommand{\hsp}{$P \cup_{\Sigma} Q$}
\newcommand{\hhs}{$V' \cup_{S'} W'$}
\newcommand{\gi}{$\gamma_{A_i}$}
\newcommand{\g}{$\gamma$}
\newcommand{\Ai}{$A_i$}
\newcommand{\scA}{${\mathcal{A}}$}
\newcommand{\AD}{${\mathcal{A}} \cap \Delta$}
\newcommand{\mA}{\mbox{\scA}}

\psfrag{V}{$V$}
\psfrag{W}{$W$}
\psfrag{Vp}{$V'$}
\psfrag{Wp}{$W'$}
\psfrag{alpha}{$\alpha$}
\psfrag{S}{$S$}
\psfrag{F}{$F$}
\psfrag{N}{$N(D_i)$}
\psfrag{L}{{\large $D$}}
\psfrag{gamma}{$\gamma$}
\psfrag{Dp}{$D'$}
\psfrag{Dpp}{$D''$}
\psfrag{b}{{\color{blue}$\beta$}}
\psfrag{Dl}{$D^{\ell}(\gamma)$}
\psfrag{A}{$A_j$}
\psfrag{gammaj}{$\gamma_j$}
\psfrag{Dj}{$D_j$}
\psfrag{N1}{$N_0$}
\psfrag{N2}{$N_1$}
\psfrag{ap}{$\alpha'$}
\psfrag{hs}{\hs}
\psfrag{hhs}{\hhs}
\psfrag{C}{$\cong$}
\psfrag{E}{$D_i(\gamma_i)'$}
\psfrag{tau}{$\tau_i$}
\psfrag{lambda}{$\gamma_i$}
\psfrag{T}{$T$}
\psfrag{tau1}{$\tau_1$}
\psfrag{E1}{$D_1(\gamma_1)'$}
\psfrag{lambda1I}{{\tiny $\gamma_1 \times \{1\}$}}
\psfrag{taui}{$\tau_i$}
\psfrag{taus}{$\tau_s$}
\psfrag{E2}{{\tiny $D_s(\gamma_s)'$}}
\psfrag{tauk}{$\tau_k$}
\psfrag{K}{$D_k(\gamma_k)'$}

\begin{abstract}
Let $T$ be a separating incompressible torus in a 3-manifold $M$. Assuming that a genus $g$ Heegaard splitting \hs \ can be positioned nicely with respect to $T$ ({\em e.g.}~ \hs \ is strongly irreducible), we obtain an upper bound on the number of stabilizations required for \hs \ to become isotopic to a Heegaard splitting which is an amalgamation along $T$. In particular, if $T$ is a canonical torus in the JSJ decomposition of $M$, then the number of necessary stabilizations is at most $4g-4$. As a corollary, this establishes an upper bound on the number of stabilizations required for \hs \ and any Heegaard splitting obtained by a Dehn twist of \hs \ along $T$ to become isotopic.
\end{abstract}

\section{Introduction}
Recent study of Heegaard splittings indicates that generically, Heegaard splittings of Haken manifolds are amalgamations along incompressible surfaces; that is, they can be decomposed into Heegaard splittings of the manifolds obtained by cutting along those surfaces (see {\em e.g.}~ \cite{BachSchSedg} for the genus $1$ case, and \cite{BachmanSchleimer}, \cite{Lackenby}, \cite{Lackenby2} and \cite{Souto} for the genus $\geq 2$ case). There are many examples, however, of Heegaard splittings (such as strongly irreducible splittings) that are not of this nature. In this paper we investigate the question of how many stabilizations are needed to make a Heegaard splitting isotopic to an amalgamation along an incompressible torus.

The peculiarity of Heegaard splittings of 3-manifolds containing incompressible tori can be seen in the recent establishment of the generalized Waldhausen conjecture by Jaco and Rubinstein \cite{JacoRubinstein} and also Li \cite{Li}, which states that a 3-manifold has only finitely many Heegaard splittings of a given genus up to isotopy, assuming the 3-manifold contains no incompressible tori. If, however, a 3-manifold contains an incompressible torus $T$, then taking a given Heegaard splitting and Dehn twisting along $T$ can yield infinitely many Heegaard splittings of the same genus (see {\em e.g.}~\cite{BDT}). Hence, a question one may ask is how many stabilizations are needed for these Heegaard splittings to become isotopic. 

Upon consideration of these questions, we have the following results (see below for relevant definitions).

\begin{theorem}
\label{MainTheorem}
Let $M$ be a closed, orientable 3-manifold, and let $T$ be a separating incompressible torus in $M$. If \hs \ is a Heegaard splitting of $M$ that can be isotoped so that $V \cap T$ consists of $k$ annuli, then after at most $k$ stabilizations, \hs \ is isotopic to an amalgamation along $T$. 
\end{theorem}

In particular, if \hs \ is strongly irreducible then (after possible isotopy) the hypotheses of the theorem are satisfied as $S$ can always be isotoped to intersect $T$ in essential simple closed curves (see {\em{e.g.}}~\cite{Schultens}). We can further refine the above bound by restricting our choice of $T$. 

\begin{corollary}
\label{GenusBound}
Let $M$ be a closed, orientable, irreducible 3-manifold, and let $T$ be a separating canonical torus in the JSJ decomposition of $M$. If \hs \ is a genus $g$ Heegaard splitting of $M$ that can be isotoped so that $V \cap T$ consists of annuli, then after at most $4g -4$ stabilizations \hs \ is isotopic to an amalgamation along $T$. 
\end{corollary}

The assumption that $M$ is irreducible is added here as it is required by the definition of the JSJ decomposition of $M$ (see Section~\ref{countingannuli}). As an immediate consequence of the above results we obtain: 

\begin{corollary}
\label{DehnTwist}
Let $M$ be a closed, orientable (irreducible) 3-manifold, and let $T$ be a separating incompressible torus (canonical torus in the JSJ decomposition) in $M$. Suppose $V \cup_S W$ and $P \cup_{\Sigma} Q$ are genus $g$ Heegaard splittings of $M$ such that $P \cup_{\Sigma} Q$ is obtained from $V \cup_S W$ via any power of a Dehn twist along $T$. Moreover, assume that \hs \ can be isotoped so that $V \cap T$ consists of $k$ annuli. Then $V \cup_S W$ and $P \cup_{\Sigma} Q$ are isotopic after at most $k$ stabilizations ($4g -4$ stabilizations). 
\end{corollary}

\section{Definitions}

Let $M$ be a closed, orientable 3-manifold. Definitions of standard terms regarding 3-manifolds can be found for example in \cite{Jaco} and \cite{Rolfsen}.

\begin{definition}
A {\em Heegaard splitting} for $M$ is a decomposition of $M$ into two handlebodies $V$ and $W$ of the same genus such that $M$ is obtained as the identification space of $V$ and $W$ identified along their boundaries via some homeomorphism from $\partial V$ to $\partial W$.
\end{definition}

\noindent The closed orientable surface $S = \partial V = \partial W$ is called the {\em splitting surface}, and we write this Heegaard splitting as $V \cup_S W$. For convenience, we do not distinguish between $V$, $W$ and $S$, and their respective embeddings in $M$. The {\em genus} of \hs \ is defined to be the genus of $S$. Two Heegaard splittings $V \cup_S W$ and $P \cup_{\Sigma} Q$ are said to be {\em isotopic} if there exists an isotopy of $M$ taking $V$ to $P$.

Given a Heegaard splitting $V \cup_S W$ of $M$, one can generate new Heegaard splittings of $M$. Let $\alpha$ be a properly embedded, boundary parallel arc in one of the handlebodies, say $V$.  Create a new handlebody $W'$ in $M$ by attaching a 1-handle $X$ to $W$ along the boundary such that $\alpha$ is the core of $X$. As $\alpha$ is boundary parallel in $V$, $V' = \overline{V - X}$ is also a handlebody. 

\begin{definition}
\rm The Heegaard splitting \hhs \ resulting from the above process is called a {\em stabilization} of \hs . A Heegaard splitting which is obtained by a stabilization of another Heegaard splitting is called {\em{stabilized}}.\end{definition}

\begin{figure}[h]
\centering
\includegraphics[width=4in]{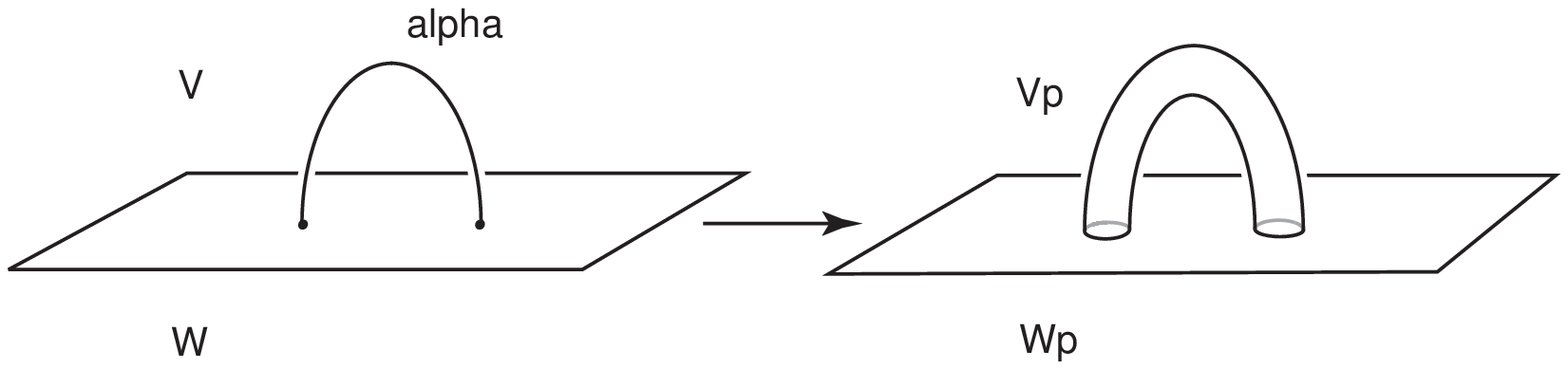}
\caption{A stabilization of \hs .} 
\label{fig:stabilization}
\end{figure}

\noindent Note that $V' \cup_{S'} W'$ has genus one larger than \hs . Repeating the above process $k$ times yields a {\em $k$-times stabilization} of \hs , or \hs \ {\em stabilized $k$ times}. It is a nice exercise to show that a Heegaard splitting being stabilized is equivalent to the property that there exist properly embedded essential disks in each of the handlebodies that intersect in a single point. Moreover, two Heegaard splittings obtained by stabilization of the same splitting are isotopic. A classical theorem of Reidemeister \cite{Reidemeister} and Singer \cite{Singer} states that any two Heegaard splittings of the same manifold can each be stabilized some indefinite number of times to become isotopic. 

\begin{definition}
\rm A Heegaard splitting \hs \ is said to be {\em reducible} if there exist essential disks $D \subset V$ and $E \subset W$ such that $\partial D = \partial E$. Otherwise $V \cup_S W$ is called {\em irreducible}. 
\end{definition} 

\begin{definition}
\rm A Heegaard splitting \hs \ is said to be {\em weakly reducible} if there are essential disks $D \subset V$ and $E \subset W$ such that $\partial D \cap \partial E = \emptyset$. Otherwise \hs \ is called {\em strongly irreducible}.
\end{definition}

\noindent A reducible Heegaard splitting is easily seen to be weakly reducible. As referred to in the introduction, a useful property of a strongly irreducible Heegaard splitting is that it can be isotoped to intersect an incompressible surface in simple closed curves essential on both surfaces. A result of Haken \cite{Haken} implies that irreducible Heegaard splittings arise only in irreducible manifolds. 

\begin{definition}
\rm Let $F$ be a separating incompressible surface in $M$. A Heegaard splitting \hs \ of $M$ is called an {\em amalgamation along $F$} if (after isotopy) $F$ is obtained from simultaneous compressions on both sides of $S$. 
\end{definition}

\begin{figure}[h]
\centering
\includegraphics[width=5in]{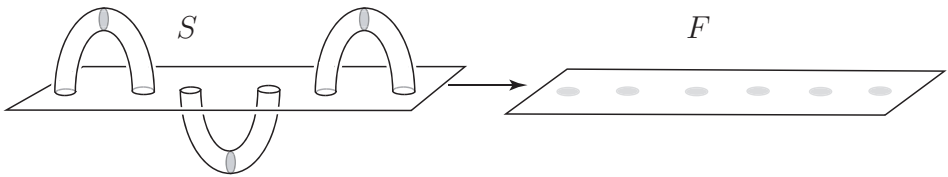}
\caption{An amalgamation along $F$.} 
\label{fig:amalgamation}
\end{figure}

\noindent This definition of amalgamation only makes sense in the context that $M$ is closed. The more general definition can be found {\em e.g.}~in \cite{Schultens2}. 

The condition that $F$ is obtained from simultaneous compressions on both sides of $S$ is equivalent to saying that $S$ can be obtained from $F$ by a series of ambient 1-surgeries on pairwise disjoint arcs properly embedded in $M$ cut along $F$. Also note that amalgamation is not unique; \hs \ can be an amalagamation along several different surfaces in $M$. By Casson and Gordon \cite{CG}, an irreducible Heegaard splitting is either strongly irreducible or an amalgamation along some incompressible surface.

\section{Amalgamations}

The purpose of this section is to determine when a Heegaard splitting is an amalgamation along a given incompressible torus. Let $A$ be an annulus properly embedded in a handlebody $V$. We say $A$ is {\em essential} in $V$ if it is incompressible and not boundary parallel in $V$. A {\em spanning arc} of $A$ is an essential arc in $A$, {\em i.e.}~a properly embedded arc in $A$ that cuts $A$ into a disk. A {\em spanning disk} for $A$ in $V$ is a boundary compressing disk for $A$ in $V$, {\em i.e.}~a disk $D$ such that $\partial D = a \cup b$ where $a = D \cap A$ is a spanning arc of $A$ and $b = D \cap \partial V$. Note that this implies $D$ is essential in $V$ cut along $A$. 

\begin{lemma} 
\label{spanningdiskexistencelemma}
Any essential annulus $A$ in a handlebody $V$ has a spanning disk. 
\end{lemma}

\begin{proof}
Using Van Kampen's Theorem one can compute a presentation for $\pi_1 (V)$ using the components of $V$ cut along $A$. As $\pi_1(V)$ is a free group, the generator of $\pi_1(A)$ must be mapped to a generator in the fundamental group of one of the components of $V$ cut along $A$. This implies that there exists an essential disk in this component which meets $A$ in a single spanning arc. Such a disk is a spanning disk of $A$ in $V$.
\end{proof}

\begin{lemma}[The Amalgamation Lemma]
\label{TheAmalgamationLemma}
Suppose $T$ is a separating incompressible torus in $M$ so that $M = N_0 \cup_T N_1$, and let \hs \ be a Heegaard splitting of $M$. 
Then $V \cup_S W$ is an amalgamation along $T$ if and only if there is an isotopy of $S$ to a surface intersecting $T$ in essential simple closed curves such that each (annulus) component of $V \cap T$ has a spanning disk in $V$ contained in $N_{\varepsilon}$, and each (annulus) component of $W \cap T$ has a spanning disk in $W$ contained in $N_{\varepsilon'}$ where $\{ \varepsilon, \varepsilon' \} = \{0,1 \}$.
\end{lemma} 

\begin{figure}[h]
\centering
\includegraphics[width=2.5in]{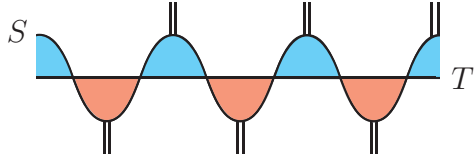}
\caption{A schematic for the Amalgamation Lemma.} 
\label{fig:amalgamation1}
\end{figure}

\begin{proof}
Assume first that \hs \ is an amalgamation along $T$. Then $S$ can be obtained from $T$ by ambient 1-surgery on pairwise disjoint arcs properly embedded in $N_0$ and $N_1$. These arcs can be isotoped by arc slides so that on each side of $T$ only one arc meets $T$ in one of its ends (see Figure~\ref{fig:tubeside}). Thus after isotopy, $S \cap T = \overline{T - (D_1 \cup D_2)}$, where $D_1$ and $D_2$ are disks in $T$.
\begin{figure}[h] 
   \centering
   \includegraphics[width=4.5in]{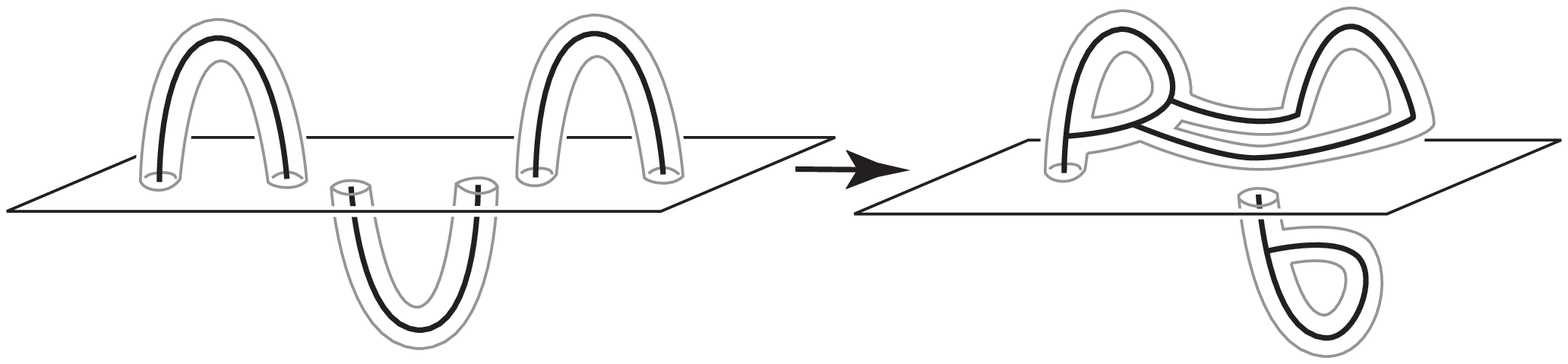} 
   \caption{\ }
   \label{fig:tubeside}
\end{figure}
Isotope $S$ into $N_0$ and $N_1$ as in Figure~\ref{fig:SecondAmalgamation} so that $S \cap T$ consists of 2 simple closed curves essential in $T$. They are essential in $S$ as well by the fact that $T$ is incompressible. The existence of spanning disks on opposite sides of $T$ for the resulting 2 annuli is then obvious, as in Figure~\ref{fig:SecondAmalgamation}.
\begin{figure}[h]
\centering
\includegraphics[width=5in]{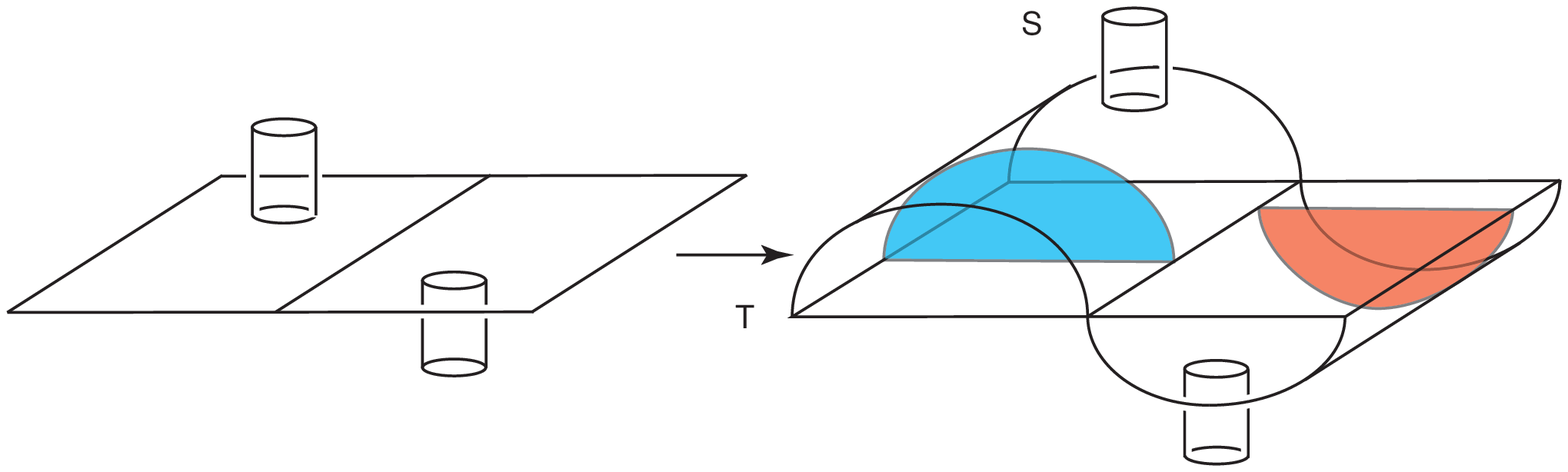}
\caption{}
\label{fig:SecondAmalgamation}
\end{figure}

For the other direction, assume without loss of generality that each component of $V \cap T$ has a spanning disk in $V$ contained in $N_0$, and each component of $W \cap T$ has a spanning disk in $W$ contained in $N_1$. Let $A_1, \ldots, A_n$ be the components of $V \cap T$ with respective spanning disks $D_1, \ldots, D_n$, and $B_1, \ldots, B_n$ the components of $W \cap T$ with respective spanning disks $E_1, \ldots, E_n$. By the definition of spanning disk and the fact that $D_i$ is contained in $N_0$, it follows that $D_i \cap T$ is a single spanning arc of $A_i$, for $1\leq i \leq n$. Similarly, $E_i \cap T$ is a single spanning arc of $B_i$, for $1 \leq i \leq n$. 

By a standard innermost disk, outermost arc argument, we may assume that $D_1, \ldots, D_n, E_1, \ldots, E_n$ are all pairwise disjoint. Each spanning disk $D_i$ defines an isotopy of $S$ in the following way. Let $N(D_i)$ be a neighborhood of $D_i$ such that $N(D_i) \cap S = \partial N(D_i) \cap S$ is a neighborhood of $\partial D_i \cap S$ in $S$, and $N(D_i) \cap T  = \partial N(D_i) \cap T$ is a neighborhood of $\partial D_i \cap T$ in $T$. Then $S$ can be isotoped so that $N(D_i) \cap S$ is replaced by $\partial N(D_i) - (N(D_i) \cap S)$. After further isotopy, $S$ intersects $T$ in a punctured annulus (see Figure~\ref{fig:ThirdAmalgamation}). After performing such an isotopy for each of the spanning disks $D_1, \ldots, D_n, E_1, \ldots, E_n$, $S$ is such that $T - S$ is a disjoint union of open disks.

\begin{figure}[h]
\centering
\includegraphics[width=5in]{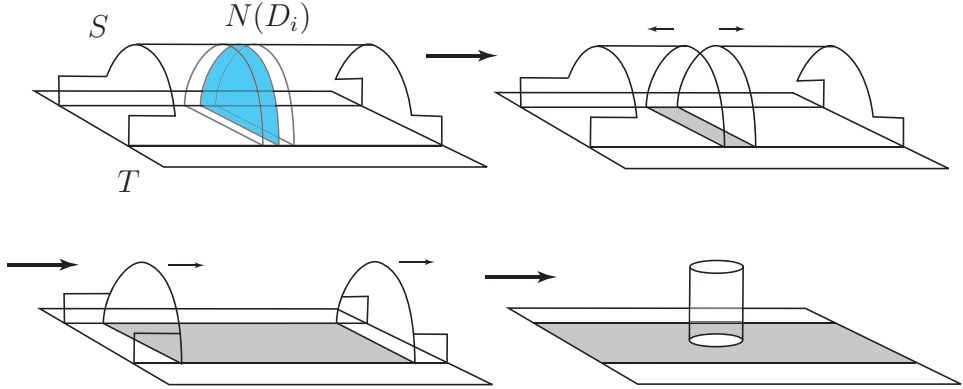}
\caption{Isotoping $S$ using the disk $D_i$.}
\label{fig:ThirdAmalgamation}
\end{figure}

If any such open disk $U$ is such that $\partial U$ bounds a disk $U_S$ in $S - T$, then by the fact that a handlebody is irreducible, $U \cup U_S$ bounds a ball. Isotope $U_S$ through this ball so that its interior equals $U$. Now, for a remaining component $V'$ of $V \cap N_0$, $V' \cap T$ is a disjoint union of disks whose interiors are a subset of the disks in $T - S$. As $V'$ is a component of $V$ cut along disks, $V'$ is a handlebody and hence has a spine. After arc slides, each arc of this spine has its endpoints on $T$. Doing this for each component of $V \cap N_0$ and $W \cap N_1$ shows that $S$ is obtained from $T$ via ambient 1-surgery along these slid arcs of the spines, implying that \hs \ is an amalgamation along $T$.
\end{proof}

\begin{remark}
\label{2kremark}
Suppose that a Heegaard splitting \hs \ can be isotoped so that $V \cap T$ consists of a disjoint union of $k$ annuli. Lemma~\ref{TheAmalgamationLemma} indicates that a $2k$-times stabilization of \hs \ is an amalgamation along $T$ in the following way. Add a 1-handle to either $V$ or $W$ such that its core is a spanning arc of a component of $W \cap T$ or $V \cap T$, respectively. The cocore of the 1-handle and a spanning disk for the annulus (which exists by Lemma~\ref{spanningdiskexistencelemma}) intersect in a single point, showing that the resulting splitting is a stabilization of \hs . By adding 1-handles disjoint from each other in this manner for every component of $T$ cut along $S$ (of which there are $2k$) and then pushing the handles added to $V$ into $N_1$ and the handles added to $W$ into $N_0$, we get the spanning disks needed for Lemma~\ref{TheAmalgamationLemma} to apply. 
\end{remark}

The purpose of the next sections is to show that the $k$ 1-handles added to one of the handlebodies in Remark~\ref{2kremark} are in fact unnecessary to obtain an amalgamation. 

\section{Annuli in handlebodies}

This section provides the necessary technical arguments used in the proof of Theorem~\ref{MainTheorem} to find spanning disks for annuli in a handlebody. Suppose a Heegaard splitting \hs \ of $M$ can be isotoped so that $V \cap T$ is a disjoint union \scA \ of annuli. Then by removing any boundary parallel components via additional isotopy we can assume the components of \scA \ are essential in $V$. Let $\Delta$ be a complete system of meridian disks for $V$. 

\begin{remark}
\label{RedefineDelta}
\rm Applying a standard cut and paste argument to $\Delta$ if necessary, we can assume that each annulus component of \scA \ meets $\Delta$ nontrivially, and so that ${\mathcal{A}} \cap \Delta$ consists of arcs properly embedded in $\Delta$, each arc being a spanning arc of some annulus in \scA. 
\end{remark}

For $D$ a component of $\Delta$, $D \cap {\mathcal{A}}$ is a collection of properly embedded arcs $\alpha_1, ... , \alpha_r$ in $D$.  

\begin{definition}
\label{Level}  
\rm Suppose that \g \ is an arc in $D \cap {\mathcal{A}}$. Then \g \ divides $D$ into two disks $D'$ and $D''$. Let $\mathcal{B}'$ be the set of properly embedded arcs in $D'$ with one endpoint on \g \ and the other endpoint on $\partial D$, intersecting each $\alpha_i$ in at most one point. Define

$$\ell_{D'} (\gamma) = \max_{\beta \in \mathcal{B}'} \{ | \beta \cap (\cup_{i=1}^r \alpha_i) |  \}.$$

\noindent Define $\ell_{D''}($\g$)$ similarly. Then we define the {\em level of $\gamma$ (in $D$)}\ to be $$\ell(\gamma) = \min \{\ell_{D'}(\gamma), \ell_{D''}(\gamma)\}.$$

For $A$ a component of $\mA$, define the {\em level of $A$ (with respect to $\Delta$)} to be $$\ell (A) = \min \{ \ell (\gamma) \}$$ where $\gamma$ is an arc component of  $A \cap \Delta$. 
\end{definition}

\noindent Note that an outermost arc for a disk $D$ is a level one arc. 

\begin{figure}
\centering
\includegraphics[width=2.3in]{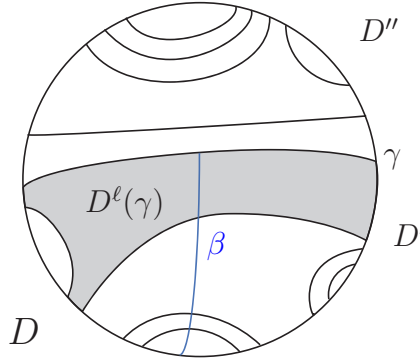}
\caption{The level adjacent disk component and an arc $\beta \in {\mathcal{B}}'$ realizing $\gamma$ as a level 4 arc in $D$.} 
\label{fig:level}
\end{figure}

\begin{definition}
\rm Let $A$ be a component of \scA \ and let $D$ be a spanning disk for $A$ such that $D \cap A = \gamma$. Define the {\em adjacent disk component of $\gamma$ with respect to $D$}, denoted $D(\gamma)$, to be the component of $D$ cut along \scA \ that contains $\gamma$. 

Let \g \ be an arc in \AD, and suppose $D$ is the disk component of $\Delta$ containing \g. Define the {\em level adjacent disk component of $\gamma$}, denoted $D^{\ell}(\gamma)$, to be the adjacent disk component of $D'$ or $D''$, whichever one is where \g \ realizes its level (if they both realize the level, choose one). See Figure~\ref{fig:level}.
\end{definition}

\indent Given a separating incompressible torus $T$ in $M$, $T$ divides $M$ into two components $N_0$ and $N_1$, where $\partial N_0 = \partial N_1 = T$. Since $V \cap T$ is assumed to be a disjoint union of essential annuli, and since $T$ is separating in $M$, each component of $V \cap T$ has an arbitrarily small neighborhood in $V$ intersecting each of $N_0$ and $N_1$ nontrivially. This gives rise to the following definition.

\begin{definition} 
\label{MutuallySeparating}
\rm Let ${\mathcal{A}}$ be a disjoint union of annuli properly embedded in a handlebody $V$. We shall say that ${\mathcal{A}}$ is {\em mutually separating} if $V$ cut along \scA \ consists of two (possibly disconnected) 3-manifolds $N_0$ and $N_1$ such that each annulus in \scA \ has an arbitrarily small neighborhood in $V$ intersecting both $N_0$ and $N_1$ nontrivially. 
\end{definition}

\begin{remark}	
\rm The condition that ${\mathcal{A}}$ is mutually separating in $V$ is equivalent to the statement that $({\mathcal{A}}, \partial {\mathcal{A}})$ represents the trivial element in $H_2(V,\partial V; \mathbb{Z}/2 \mathbb{Z})$.
\end{remark}

\noindent Note that \scA \ can be mutually separating even if some of its component annuli are themselves nonseparating in $V$.

\begin{lemma}
\label{ordering}
Let ${\mathcal{A}}$ be a mutually separating disjoint union of essential annuli properly embedded in a handlebody $V$. Then there is an ordering $A_1, \ldots, A_{k}$ of the annuli in ${\mathcal{A}}$ such that for each $1 \leq i \leq k$, $D_{i}$ is a spanning disk and $\gamma_i$ is a spanning arc for $A_i$, and

$$D_i \cap A_j = \left\{ \begin{array}{cl}
	\mbox{\parbox[t]{2in}{(a possibly empty set of) arcs properly embedded in $D_i$ and parallel  in $ A_j$ to $\gamma_j$}} & \mbox{if $j < i$} \\
	\ & \ \\
	\gamma_i \subset \partial D_i& \mbox{if $j=i$} \\
	\ & \ \\
	\emptyset & \mbox{otherwise} 
	\end{array}	\right. $$

Furthermore, suppose $\alpha$ ($\neq \gamma_i$) is an arc in $\mbox{\scA} \cap D_i$ contained in $\partial D_i(\gamma_i)$, so $\alpha$ is parallel in $A_j$ to $\gamma_j$ for some $j < i$. If $D_i(\gamma_i) \subset N_{\varepsilon}$, then $D_j(\gamma_j) \subset N_{\varepsilon'}$ where $\{ \varepsilon, \varepsilon' \} = \{0,1 \}$. 
\end{lemma}

\begin{proof} 
Let $\Delta$ be a complete system of meridian disks for $V$, chosen as in Remark~\ref{RedefineDelta} so that \AD \ is transverse and consists of spanning arcs of the annuli. We will order the annuli in \scA \ and construct the disks $D_i$ using the level of the annuli with respect to $\Delta$. 

First suppose $A$ is a component of \scA \ which is level one. That is, there is a level one arc \g \ in $A \cap \Delta$. Label $A$ as $A_1$, $\gamma$ as $\gamma_1$, and take $D_1$ to be the level adjacent disk component of $\gamma$ ({\em i.e.}~the outermost disk cut off by $\gamma$ of the component of $\Delta$ containing $\gamma$). Continue this for all level one annuli which are components of \scA. We obtain a list of annuli $A_1, \ldots ,A_{k_1}$ with respective spanning arcs $\gamma_1, \ldots , \gamma_{k_1}$, and respective disks $D_1, \ldots ,D_{k_1}$. Note that each $D_i$ for $1\leq i \leq k_1$ is contained in $N_0$ or $N_1$, thereby satisfying the conclusion of the lemma. 

Continuing in an inductive manner, suppose that we have compiled a list $A_{1}, \ldots, A_{k_{m-1}}$ of annuli of level $\leq m-1$, along with their corresponding spanning arcs $\gamma_{1}, \ldots, \gamma_{k_{m-1}}$ and disks $D_{1}, \ldots, D_{k_{m-1}}$. Proceeding as before, let $A$ be a level $m$ annulus in \scA, and let $\gamma$ be a level $m$ arc in $A \cap \Delta$. Set $A = A_{k_{m-1}+1}$, and set $\gamma = \gamma_{k_{m-1}+1}$. It remains to construct $D_{k_{m-1} + 1}$.
 
Let $D$ be the component of $\Delta$ containing $\gamma$. Assume without loss of generality that $D^{\ell}(\gamma)$ is contained in $N_0$. Let $\alpha$ be an arc in $\partial D^{\ell}(\gamma) \cap \mA$ other than $\gamma$. Then $\ell (\alpha)< m$ and hence $\alpha$ is a spanning arc of some annulus $A_j$, $j  \leq k_{m-1}$. The annulus $A_j$ has corresponding spanning arc $\gamma_j$. First assume that $\alpha$ is the only such arc other than $\gamma$ contained in $\partial D^{\ell}(\gamma) \cap \mA$. 

Assume first that $D_j(\gamma_j)$ is contained in $N_0$. Since $\alpha$ and $\gamma_j$ are spanning arcs of $A_j$, they cut off a rectangle in $A_j$ (if $\alpha = \gamma_j$, then this rectangle is simply the arc $\alpha$). Taking that rectangle and attaching $D_j$ and $D^{\ell}(\gamma)$ gives a disk, which can be pushed slightly into $N_0$. Take this disk as $D_{k_{m-1} + 1}$. Note that $\alpha$ has been eliminated as an arc of intersection in $D_{k_{m-1} + 1}$ (see Figure~\ref{fig:adjdisk1}). Since the disk $D_j$ satisfies the conclusion of the lemma, this implies $D_{k_{m-1} + 1}$ does as well.

\begin{figure}[h]
\centering
\includegraphics[width=4in]{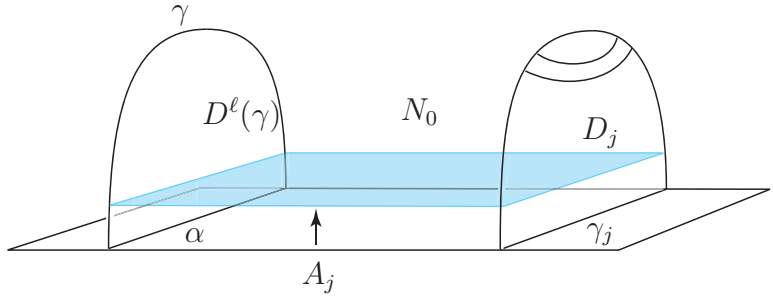}
\caption{$D_j(\gamma_j)$ is in $N_0$.} 
\label{fig:adjdisk1}
\end{figure}

Now assume that $D_j(\gamma_j)$ is contained in $N_1$. If $\alpha = \gamma_j$, set $D_{k_{m-1} + 1} = D^{\ell}(\gamma) \cup D_j$. If $\alpha \neq \gamma_j$, as above take a rectangle in $A_j$ bounded by $\alpha$ and $\gamma_j$. The rectangle, along with $D^{\ell}(\gamma)$ and $D_j$, forms a disk. Isotope the disk by pushing the subdisk formed by $D^{\ell}(\gamma)$ and the rectangle slightly into $N_0$, keeping $D_j$ fixed. Take the resulting disk as $D_{k_{m-1} + 1}$ (see Figure~\ref{fig:adjdisk2}). As before, $D_{k_{m-1} + 1}$ satisfies the conclusion of the lemma since $D_j$ does. 

\begin{figure}[h]
\centering
\includegraphics[width=4in]{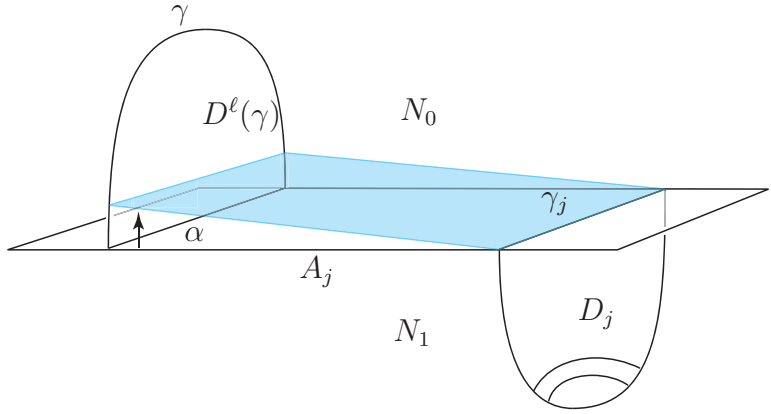}
\caption{$D_j(\gamma_j)$ is in $N_1$.} \label{fig:adjdisk2}
\end{figure}

Now suppose $\alpha'$ is an arc in addition to \g \ and $\alpha$ in $\partial D^{\ell}(\gamma) \cap \mA$. If $\alpha'$ is not a spanning arc of the same annulus component of $\mathcal{A}$ as $\alpha$, then the above construction applies similarly to $\alpha'$ to obtain $D_{k_{m-1}+1}$. If, however, $\alpha'$ is a spanning arc of the same annulus as $\alpha$, then a slight modification of the above argument is needed in order to attach a disk to $\partial D^{\ell}(\gamma)$ at $\alpha'$. Let $A_j$ be the annulus containing $\alpha$ and $\alpha'$, and, as above, consider the disk $D_j$. As before, attach $D_j$ to $D^{\ell}(\gamma)$ along a rectangle between $\alpha$ and $\gamma_j$, and then isotope appropriately off of $A_j$. Then for $\alpha'$, take a parallel copy $D_j'$ of $D_j$ in $V$ cut along $A_j$, and a rectangle from $\alpha'$ to the arc in $\partial D_j'$ parallel to $\gamma_j$, and push off into $N_0$ as before, depending on whether $D_j(\gamma_j)$ is in $N_0$ or $N_1$. By choosing $D_j'$ to be on the appropriate side of $D_j$, we can ensure that the rectangle between $\alpha'$ and $D_j'$ is disjoint from the rectangle between $\alpha$ and $D_j$. Repeating this process for any additional arcs in $\partial D^{\ell}(\gamma) \cap \mA$ (except for \g), we obtain the desired disk $D_{k_{m-1}+1}$ (see Figure~\ref{fig:adjdisk3}).

\begin{figure}[h]
\centering
\includegraphics[width=4.5in]{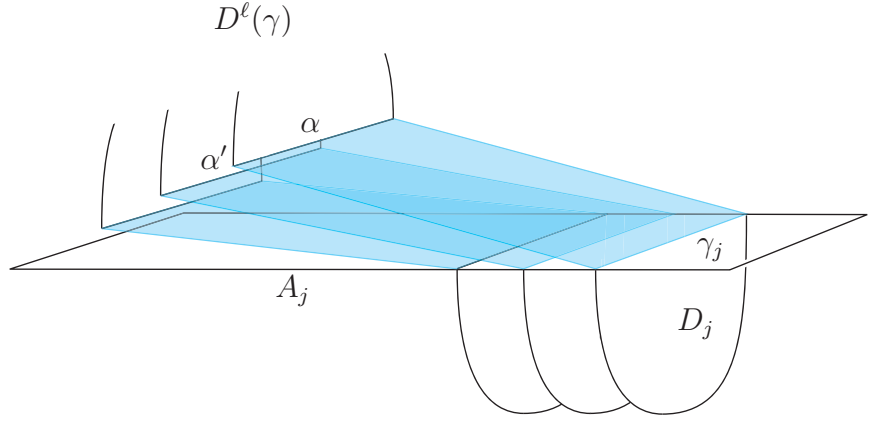}
\caption{$\partial D^{\ell}(\gamma) \cap \mA$ contains more than one arc other than \g.} 
\label{fig:adjdisk3}
\end{figure}

\end{proof}

\section{Proof of the main theorem}

\subsection{Proof of Theorem~\ref{MainTheorem}}

Note that both $W \cap T$  and $V \cap T$ are nonempty since $T$ is incompressible and cannot be contained in a handlebody.  Let $\lambda_1, \ldots, \lambda_k$ be spanning arcs of the components of $W \cap T$. Attach $k$ 1-handles to $V$ so that their cores are $\lambda_1, \ldots, \lambda_k$, respectively. As discussed in Remark~\ref{2kremark} this yields a $k$-times stabilization $V' \cup_{S'} W'$ of \hs.

Let $T \times I$ be a regular neighborhood of $T$ such that $T \times \{ 1/2 \} = T$. Then \hs \ and \hhs \ can be considered identical except for inside $T \times I$. 

\begin{figure}[h]
\centering
\includegraphics[width=5in]{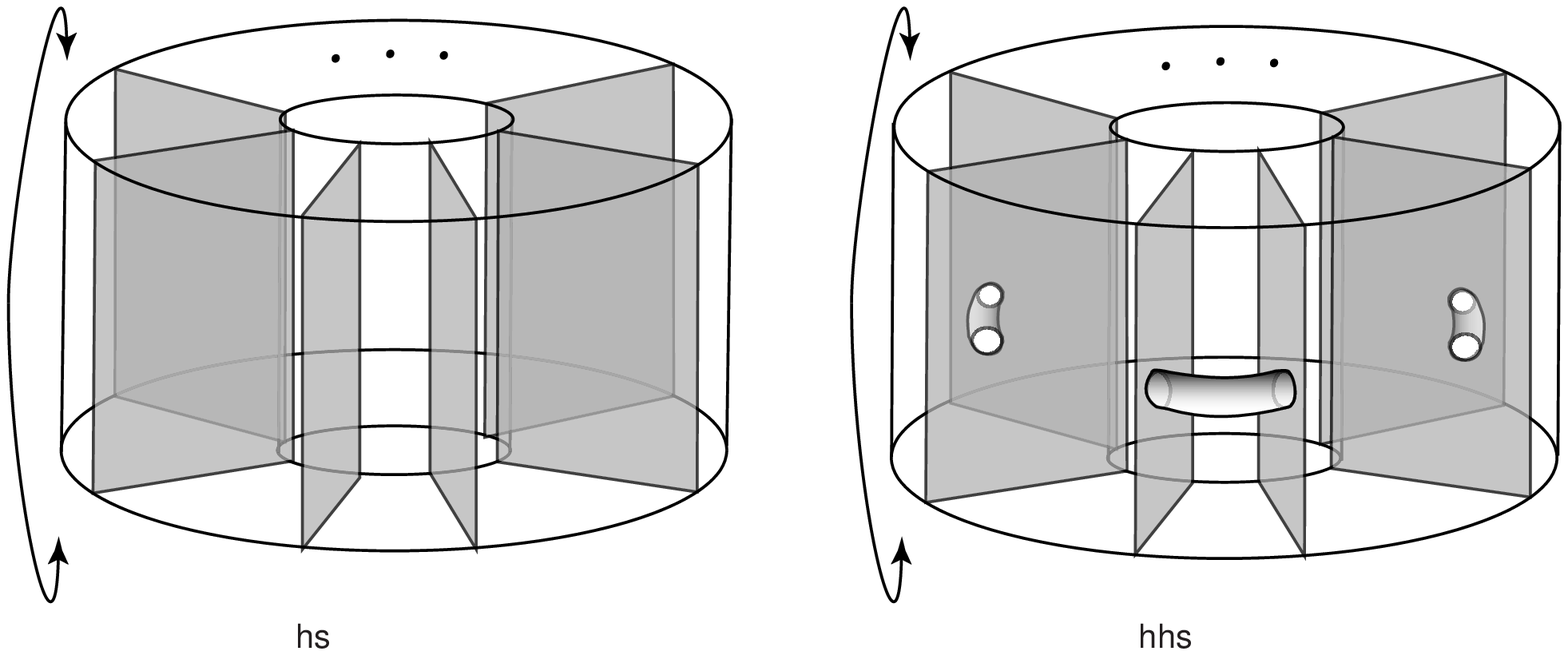}
\caption{\hs \ and \hhs \ inside $T \times I$.} 
\label{fig:stabilize}
\end{figure}

\begin{figure}[h]
\centering
\includegraphics[width=3.5in]{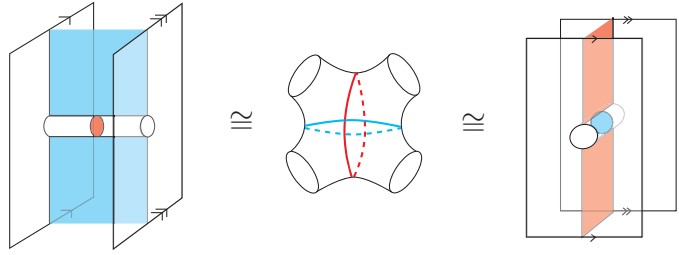}
\caption{Isotopic pictures of annuli connected by a tube.} 
\label{fig:DualAnnuli}
\end{figure}

\begin{figure}[h]
\centering
\includegraphics[width=4in]{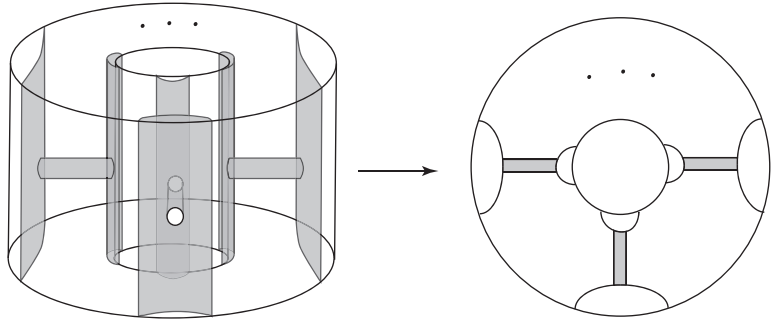}
\caption{ A schematic for $V' \cup_{S'} W'$ inside $T \times I$.} 
\label{fig:DualStabilize}
\end{figure}

A {\em tube} $\tau$ in a handlebody $V$ is a regular neighborhood $D \times I$ of a compressing disk $D$ of $V$, so that $\partial D \times I \subset \partial V$. We refer to $D \times \{0 \}$ and $D \times \{1\}$ as the {\em feet} of $\tau$. Observe that two parallel annuli connected by a tube in $M$ is isotopic to a dual picture, as in Figure~\ref{fig:DualAnnuli}. Performing this isotopy in $T \times I$ then yields a configuration of $V' \cup_{S'} W' \cap (T \times I)$ as in the first part of Figure~\ref{fig:DualStabilize}. If we take a top-down view, we obtain the schematic picture in the second part of Figure~\ref{fig:DualStabilize}. The tubes running  through $T \times I$ are the ``spokes" in this schematic.

Let \scA \ be the disjoint union of annuli in $V \cap T$. Apply Lemma~\ref{ordering} to obtain an ordering $A_1,\ldots ,A_k$ of the annuli in \scA , spanning arcs $\gamma_1, \ldots, \gamma_k$ of the annuli, and spanning disks $D_1, \dots ,D_k$ satisfying the conclusion of the lemma. Note that these disks can be chosen to miss the stabilizations added above, thus they exist for $V'$. 

\begin{figure}[h]
\centering
\includegraphics[width=3in]{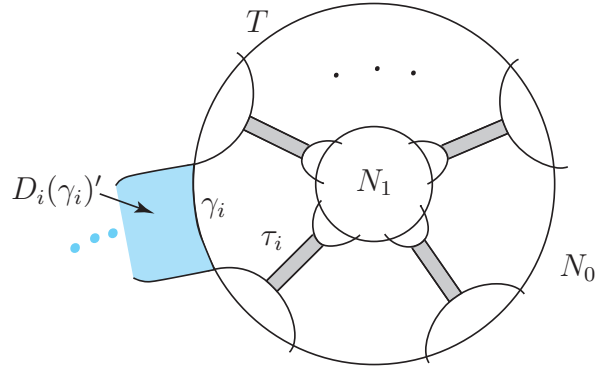}
\caption{${D_i}(\gamma_i)'$, $\tau_i$, and $\gamma_i$.} 
\label{fig:wheel2}
\end{figure}

As in the proof of Lemma~\ref{ordering}, we may assume that any arcs in $D_i \cap {\mathcal{A}}$ parallel in $A_j$ to $\gamma_j$ for $j < i$ are arbitrarily close to $\gamma_j$ in $A_j$. Choose the $D_i$ so that they intersect $T \times I$ in $\gamma_i \times I$, $1 \leq i \leq k$. For notational purposes, if ${D_i}(\gamma_i)$ is the adjacent disk component of $\gamma_i$ with respect to $D_i$, set ${D_i}(\gamma_i)' = {D_i}(\gamma_i) \cap (\overline{M - (T \times I)})$. Let $\tau_i$, $1 \leq i \leq k$, be the tubes in $W'$ resulting from the stabilizations of \hs \ as in Figure~\ref{fig:DualStabilize}, ordered so that $\tau_i$ is the tube immediately counterclockwise in the schematic from $\gamma_i$. See Figure~\ref{fig:wheel2}. 

\begin{figure}[h]
\centering
\includegraphics[width=5in]{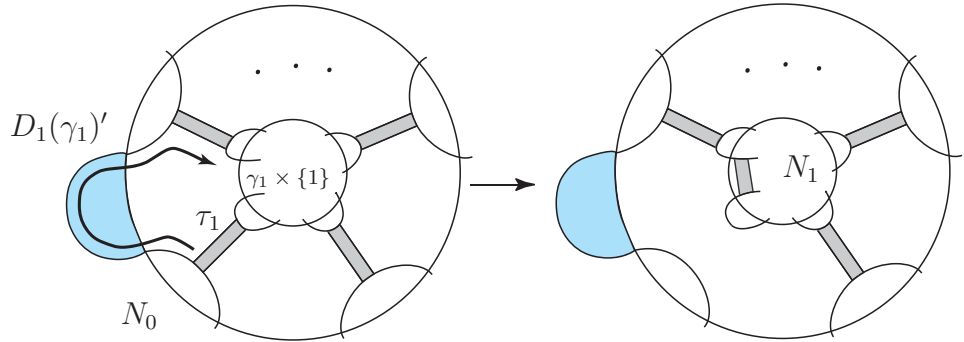}
\caption{Isotoping $\tau_1$.} 
\label{fig:wheel3}
\end{figure}

We now use the disks $D_i$ to isotope the tubes $\tau_i$ as follows. First, consider $\tau_1$ and the disk $D_1$. Note that $D_1$ is the adjacent disk component $D_1(\gamma_1)$ and lies in either $N_0$ or $N_1$. Assume first that $D_1 \subset N_0$. As $D_1$ intersects ${\mathcal{A}}$ only at $\gamma_1$, isotope the foot of $\tau_1$ lying on $T \times \{0\}$ across  a regular neighborhood of $D_1(\gamma_1)'$, and then down the tube immediately clockwise from $\tau_1$ in the schematic. The result is that $\tau_1$ is isotoped so that its core is $\gamma_1 \times \{1\}$ in $T \times I$. Push $\tau_1$ slightly into $N_1$ (we say that $\tau_1$ is {\em adjacent} to $\gamma_1 \times \{1\}$ after this isotopy). See Figure~\ref{fig:wheel3}. 

If, on the other hand, $D_1 \subset N_1$, then by a symmetric argument the other foot of $\tau_1$ can be isotoped through a neighborhood of $D_1(\gamma_1)'$ so that it is adjacent to $\gamma_1 \times \{0\}$. 

\begin{figure}[h]
\centering
\includegraphics[width=3in]{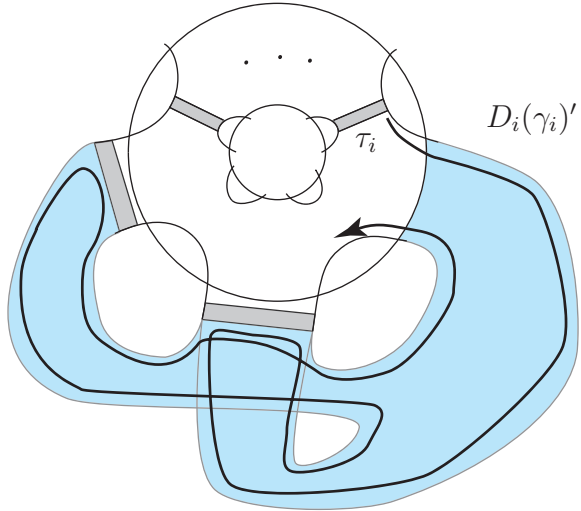}
\caption{Isotoping $\tau_i$ through $D_i(\gamma_i)'$.} 
\label{fig:wheel4}
\end{figure}

Given $2 \leq i \leq k-1$, isotope all tubes $\tau_j$, $j<i$, to lie adjacent to $\gamma_j \times \{0\}$ or $\gamma_j \times \{1\}$ in $T \times I$ as above (note that this depends on whether or not $D_j(\gamma_j)$ is in $N_1$ or $N_0$). To isotope $\tau_i$, assume that $D_i(\gamma_i)$ is in $N_0$ (as above a symmetric argument applies if $D_i(\gamma_i)$ is in $N_1$).  By Lemma~\ref{ordering}, each of the arcs $\alpha_1, ... , \alpha_m$ in $(\partial D_i(\gamma_i) \cap {\mathcal{A}}) - \gamma_i$ is parallel in $A_j$ to $\gamma_j$ for some $j < i$. Let $j_1, ... , j_n$ be this set of indices. For these $j_l$, the latter conclusion of Lemma~\ref{ordering} implies that $D_{j_l}(\gamma_{j_l})$ is contained in $N_1$. Hence, $\tau_{j_l}$ is adjacent to $\gamma_{j_l} \times \{0 \}$. Since we chose the arcs $\alpha_1, ..., \alpha_m$ to lie arbitrarily close to the spanning arc $\gamma_{j_l}$ of the annulus in which they lie, $D_i(\gamma_i)$ may be isotoped so that it meets $\partial \tau_{j_1} \cup \ldots \cup \partial \tau_{j_n}$ in place of $\alpha_1\cup\ldots \cup \alpha_m$. The upshot is that $D_i(\gamma_i)$ is a disk in $V'$ that meets $T$ only in the arc $\gamma_i$. Thus, as before, keeping one foot fixed isotope $\tau_i$ through a regular neighborhood of $D_i(\gamma_i)'$ until $\tau_i$ is in $T \times I$ and the other foot of $\tau_i$ has reached the next section clockwise in the schematic picture. See Figure~\ref{fig:wheel4}.

Having isotoped the foot of $\tau_i$ through to the next section of the schematic, one of two things can occur. If the tube $\tau_s$ originating in that section has not yet been isotoped ({\em i.e.}~$s > i$), then as before the foot of $\tau_i$ may be isotoped along $\tau_s$ causing $\tau_i$ to become adjacent to $\gamma_i \times \{ 1 \}$ as desired. If, on the other hand, $s < i$ so that $\tau_s$ has already been isotoped, then $\tau_s$ is adjacent to $\gamma_s \times \{ \varepsilon \}$ for $\varepsilon = 0$ or $1$. But this implies that $D_s(\gamma_s)'$ meets $\gamma_s \times \{1 - \varepsilon\}$ along its boundary ({\em i.e.}~$D_s(\gamma_s) \subset N_0$ if $\varepsilon =1$ or  $D_s(\gamma_s) \subset N_1$ if $\varepsilon =0$). Lemma~\ref{ordering} ensures that any arc $\alpha$ in $(\partial D_s(\gamma_s) \cap {\mathcal{A}}) - \gamma_s$ is parallel in $A_j$ to $\gamma_j$ for some $j < s$, so by the same argument as above, the corresponding components of $(\partial D_s(\gamma_s) \cap {\mathcal{A}}) - \gamma_s$ lie on already isotoped tubes. Therefore, we may continue to isotope $\tau_i$ across to the next section of the schematic, either by sliding the foot of $\tau_i$ across $\tau_s$ or through a regular neighborhood of $D_s (\gamma_s)'$. 

\begin{figure}[h]
\centering
\includegraphics[width=2.5in]{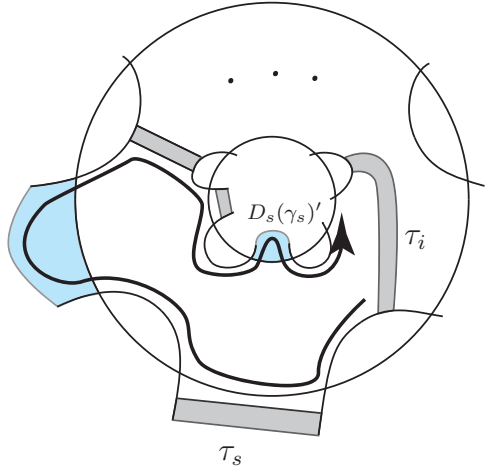}
\caption{Isotoping $\tau_i$ to be adjacent to $\gamma_i \times \{ 1 \}$.} 
\label{fig:wheel5}
\end{figure}

Continue sliding across sections in this manner until reaching a tube $\tau_{s'}$ that has not been isotoped ({\em i.e.}~$s' > i$). Slide the foot of $\tau_i$ down $\tau_{s'}$, and then back through all the previous sections, this time on the other side. That is, if the foot of $\tau_i$ was isotoped initially along $T \times \{0\}$ and across the isotoped $\tau_s$, then after having been isotoped down $\tau_{s'}$, $\tau_i$ can be isotoped along $T \times \{1\}$ and through a regular neighborhood of $D_s(\gamma_s)'$. Continue the isotopy so that $\tau_i$ becomes adjacent to $\gamma_i \times \{1\}$ as desired (see Figure~\ref{fig:wheel5}). Having done this process for each of the tubes $\tau_i$, $1 \leq i \leq k-1$, we leave the remaining tube $\tau_k$ unmoved. 

\begin{figure}[h]
\centering
\includegraphics[width=2.5in]{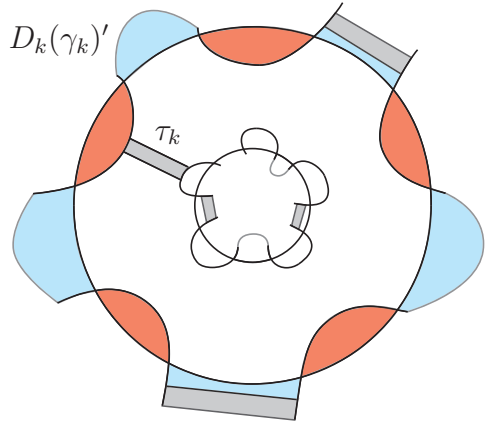}
\caption{Spanning disks for each of the annuli in $V' \cap T$ and $W' \cap T$.} 
\label{fig:wheel6}
\end{figure}

Suppose that $D_k(\gamma_k)$ is in $N_0$ (again, a symmetric argument applies if $D_k (\gamma_k)$ is in $N_1$). As above, we can take $D_k(\gamma_k)$ to be completely in $N_0$ so that it meets $T$ only at $\gamma_k$. That is, $D_k (\gamma_k)$ is a spanning disk for $A_k$ contained in $N_0$. Each of the tubes $\tau_i$, $1 \leq i \leq k-1$, is adjacent to either $\gamma_i \times \{ 0 \}$ or $\gamma_i \times \{ 1 \}$. In the former case, there is a disk between $\partial \tau_i$ and $\gamma_i \times \{ 0 \}$ which is a spanning disk for $A_i$ contained completely in $N_0$. In the latter case, as was argued above, $D_i(\gamma_i)'$ is a spanning disk for $A_i$ contained completely in $N_0$. Each of the annuli in $W \cap T$ clearly has a spanning disk contained in $N_1$, as indicated in Figure~\ref{fig:wheel6}. 
Thus, by Lemma~\ref{TheAmalgamationLemma}, \hhs \ is an amalgamation along $T$. 
\begin{flushright} 
$\Box$
\end{flushright}

\subsection{Dehn twisting and Corollary~\ref{DehnTwist}}

Consider the torus $T$ in $M$ as the product $S^1 \times S^1$, so that a point on $T$ can be written in the form $(x,y)$. Define the map $r_{\theta}: S^1 \rightarrow S^1$ to be a rotation of $S^1$ of angle $\theta$.

\begin{definition} 
\rm A {\em Dehn twist along $T$} is a homeomorphism $h \colon M \to M$ such that in a product neighborhood $T \times I$ of $T$, 
$$h (x,y,t) = (r_{2 \pi p t} (x), r_{2 \pi q t} (y), t),$$ 
where $p$ and $q$ are relatively prime integers, and $h|{\overline{M - (T \times I)}}$ is the identity. 
\end{definition}

\begin{remark}
\rm Any homeomorphism of $M$ to itself which is the identity outside a product neighborhood $T \times I$ is isotopic to a power of a Dehn twist along $T$. 
\end{remark}

\begin{proof}[Proof of Corollary~\ref{DehnTwist}] Suppose that \hs \ and \hsp \ are Heegaard splittings of $M$ such that \hs \ can be isotoped so that $V \cap T$ consists of $k$ annuli, and \hsp \ is obtained from \hs \ by any power of a Dehn twist along $T$. Theorem~\ref{MainTheorem} implies that both splittings are isotopic to amalgamations along $T$ after at most $k$ stabilizations. Let \hhs \ and $P' \cup_{\Sigma'} Q'$ be \hs \ and \hsp \ stabilized $k$ times, respectively. Let $J \colon M \times I \to M$ be the isotopy of \hhs \ constructed in the proof of Theorem~\ref{MainTheorem}, and let $h \colon M \to M$ be a power of a Dehn twist along $T$. Note that $J \circ (h \times id) = h \circ J$, which implies that $J((P' \cup_{\Sigma'} Q'), 1)$ is obtained from $J((V' \cup_{S'} W'),1)$ by a power of a Dehn twist along $T$. It is a simple exercise to observe that a power of a Dehn twist along $T$ of a Heegaard splitting which is an amalgamation along $T$ is isotopic in $T \times I$ to the original splitting. This implies that \hhs \ and $P' \cup_{\Sigma'} Q'$ are isotopic. \end{proof}

\section{Counting annuli in handlebodies}
\label{countingannuli}

In this section we adapt arguments from \cite{MS} to prove Corollary~\ref{GenusBound}. Recall a theorem of Jaco and Shalen and also Johannson (see for example \cite{Jaco}) that states in a compact, orientable, irreducible 3-manifold $M$ there exists a disjoint union of incompressible tori $\Theta$, unique up to isotopy, such that each component of $M$ cut along $\Theta$ either is a Seifert fibered space or is atoroidal. This decomposition of $M$ is called the {\em JSJ decomposition of $M$}. 

\begin{definition}
A component $T$ of $\Theta$ is called a {\em canonical torus in the JSJ decomposition of $M$}. 
\end{definition}

Any incompressible torus in $M$ is isotopic either to a canonical torus in the JSJ decomposition or to a torus contained completely in one of the Seifert fibered components. The following definitions are taken from Section 4 of ~\cite{MS}.

\begin{definition}
\rm Let ${\mathcal{A}}$ be a disjoint union of essential annuli in a handlebody $V$. A component $Z$ of $V$ cut along ${\mathcal{A}}$ is called {\em toral} if $Z$ is a solid torus. 

Define the {\em complexity} of a toral component $Z$ to be $c(Z) = | \partial V \cap Z| - \varepsilon$, where $\varepsilon = 1$ if the annuli ${\mathcal{A}} \cap Z$ are longitudes of $Z$ and $\varepsilon = 0$ otherwise. 

The complexity of a union of toral components is defined to be the sum of the complexities of the individual components.
\end{definition}
 
\begin{definition}
\rm Suppose \hs \ is a Heegaard splitting of $M$ and that $N$ is a Seifert fibered component in the JSJ decomposition of $M$. Then $N$ is called {\em aligned} with respect to \hs \ if $S$ intersects $\partial N$ only in fibers.
\end{definition}

The following theorem, due to Waldhausen, characterizes incompressible, $\partial$-incompressible surfaces in Seifert fibered spaces. The interested reader is refered to \cite{Jaco}, Theorem VI.32 for a proof.

\begin{theorem} 
\label{ja}
A properly embedded incompressible and $\partial$-incompressible two-sided surface in an orientable Seifert fibered space $N$ can be properly isotoped so that either it is vertical (a union of fibers) or it is horizontal (transverse to the fibering). If it is horizontal, then $N$ is the union of two copies of an $I$-bundle $\xi$ over a surface $E$, glued along their $\partial I$-bundles $\xi'$, and the incompressible surface consists of parallel copies of $\xi'$.
\end{theorem}

Assume as in the hypothesis of Theorem~\ref{MainTheorem} that \hs \ can be isotoped to intersect $T$ such that $V \cap T$ consists of essential annuli. We now prove the first of three results that will establish an upper bound for the number of annuli in $V \cap T$. This proof can be found in the proof of a more general result, Theorem 4.7 in ~\cite{MS}. We include it here for completeness.

\begin{lemma} 
\label{toral}
If $T$ is a separating canonical torus in the JSJ decomposition of $M$ and \hs \ is a Heegaard splitting of $M$ isotoped so that $V \cap T$ consists of a minimal number of essential annuli, then each component of $V \cap T$ intersects at most one toral component of $V$ cut along $V \cap T$.
\end{lemma}

\begin{proof} Suppose $N$ is a component of the JSJ decomposition of $M$, and that $T$ is a component of $\partial N$. Suppose $N$ contains a toral component of $V$ cut along $T$. The incompressibility of $T$ and the hypothesis that $|V \cap T|$ is minimal imply that $N$ contains an essential annulus, and hence is Seifert fibered. If $N$ is non-aligned, then the annulus cannot be vertical. By Theorem~\ref{ja}, the annulus must therefore be horizontal. Moreover, Theorem~\ref{ja} implies that $N$ is an $I$-bundle over an annulus or a M\"{o}bius band. In both cases, $N$ can be refibered to be aligned. Note that this follows since $T$ is separating, and hence does not meet $N$ on both sides. 

Hence, if two adjacent components of $V$ cut along $T$ are toral, they have to be in aligned Seifert fibered spaces. But this implies that $T$ is incident to Seifert fibered spaces that are aligned on both sides, implying $T$ is not a canonical torus in the JSJ decomposition. 
\end{proof}

\begin{lemma}
\label{complexity}
Let ${\mathcal{A}}$ be a disjoint union of essential annuli in a handlebody $V$, and let $Z$ be the union of toral components of $V$ cut along ${\mathcal{A}}$. Suppose that at most one of the components of $V$ cut along ${\mathcal{A}}$ on either side of a component of ${\mathcal A}$ is toral, and let $\alpha$ be the number of annuli which do not meet toral components (on either side). Then $c(Z) + \alpha \leq 2g - 2$, where $g$ is the genus of $\partial V$. 
\end{lemma}

This is Lemma 4.4 in \cite{MS}.

\begin{lemma}
\label{bound}
Let $T$ be a canonical torus in the JSJ decomposition of $M$ and let \hs \ be a Heegaard splitting of $M$ isotoped so that $V \cap T$ consists of a minimal number $k$ of essential annuli. Then $k \leq 4g - 4$. 
\end{lemma}

\begin{proof} 
Let $\alpha$ be the number of annuli in ${\mathcal{A}} = V \cap T$ which are not incident to toral components of $V$ cut along ${\mathcal A}$ on either side, and let $\beta$ be the number of annuli which meet a toral component. Then $k = \alpha + \beta$. Let $Z_1, ... , Z_n$ be the toral components of $V$ cut along ${\mathcal A}$, and let $\beta_i$ be the number of components of ${\mathcal{A}}$ incident to $Z_i$, $1 \leq i \leq n$. Then $c(Z_i) = \beta_i - 1$ if $Z_i \cap {\mathcal{A}}$ are longitudes, and $c(Z_i) = \beta_i$ otherwise. Since by Lemma~\ref{toral} no two adjacent components of $V$ cut along ${\mathcal A}$ are toral, we have 
$$\beta = \sum_{i=1}^n \beta_i \leq \sum_{i=1}^n (c(Z_i) + 1) = c(Z) + n,$$ where $Z = \cup_{i=1}^n Z_i$. 

Now, observe that for each toral component $Z_i$, we have $c(Z_i) \geq 1$. For if some $Z_i$ were such that $c(Z_i) = 0$, then $Z_i$ would have one longitudinal annulus on its boundary in $\partial V$, implying that the sole annulus in $Z_i \cap {\mathcal{A}}$ is boundary parallel and hence not essential. Therefore, $n \leq c(Z)$, and from the above we conclude that
$$\beta \leq 2c(Z).$$
Hence,
$$k = \alpha + \beta \leq \alpha + 2c(Z) \leq 2\alpha + 2 c(Z) \leq 4g - 4$$
by Lemma~\ref{complexity}. 
\end{proof}


The proof of Corollary~\ref{GenusBound} follows readily from Theorem~\ref{MainTheorem} and Lemma~\ref{bound}.

\bibliographystyle{plain}
\bibliography{Dehn_twist}

\begin{thebibliography}{10}

\bibitem{BDT}
David Bachman and Ryan Derby-Talbot.
\newblock Non-isotopic {H}eegaard splittings of {S}eifert fibered spaces.
\newblock {\em Algeb. Geom. Topol.}, 6:351--372, 2006.

\bibitem{BachmanSchleimer}
David Bachman and Saul Schleimer.
\newblock Surface bundles versus {H}eegaard splittings.
\newblock {\em Comm. Anal. Geom.}, 13(5):903--928, 2005.

\bibitem{BachSchSedg}
David Bachman, Saul Schleimer, and Eric Sedgwick.
\newblock Sweepouts of amalgamated 3-manifolds.
\newblock {\em Algeb. Geom. Topol.}, 6:171--194, 2006.

\bibitem{CG}
A.~J. Casson and C.~McA. Gordon.
\newblock Reducing {H}eegaard splittings.
\newblock {\em Topology Appl.}, 27(3):275--283, 1987.

\bibitem{Haken}
Wolfgang Haken.
\newblock Some results on surfaces in {$3$}-manifolds.
\newblock In {\em Studies in Modern Topology}, pages 39--98. Math. Assoc. Amer.
  (distributed by Prentice-Hall, Englewood Cliffs, N.J.), 1968.

\bibitem{Jaco}
William Jaco.
\newblock {\em Lectures on three-manifold topology}, volume~43 of {\em CBMS
  Regional Conference Series in Mathematics}.
\newblock American Mathematical Society, Providence, R.I., 1980.

\bibitem{JacoRubinstein}
William Jaco and J.~Hyam Rubinstein.
\newblock {$0$}-efficient triangulations of 3-manifolds.
\newblock {\em J. Differential Geom.}, 65(1):61--168, 2003.

\bibitem{Lackenby}
Marc Lackenby.
\newblock The {H}eegaard genus of amalgamated 3-manifolds.
\newblock {\em Geom. Dedicata}, 109:139--145, 2004.

\bibitem{Lackenby2}
Marc Lackenby.
\newblock Heegaard splittings, the virtually {H}aken conjecture and {P}roperty
  (tau).
\newblock {\em Invent. Math.}, 164:317--359, 2006.

\bibitem{Li}
Tao Li.
\newblock Heegaard surfaces and measured laminations. {I}. {T}he {W}aldhausen
  conjecture.
\newblock {\em Invent. Math.}, 167(1):135--177, 2007.

\bibitem{Reidemeister}
Kurt Reidemeister.
\newblock Zur dreidimensionalen topologie.
\newblock {\em Abh. Math. Sem. Univ. Hamburg}, 9:189--194, 1933.

\bibitem{Rolfsen}
Dale Rolfsen.
\newblock {\em Knots and links}.
\newblock Publish or Perish Inc., Berkeley, Calif., 1976.
\newblock Mathematics Lecture Series, No. 7.

\bibitem{MS}
Martin Scharlemann and Jennifer Schultens.
\newblock Comparing {H}eegaard and {JSJ} structures of orientable 3-manifolds.
\newblock {\em Trans. Amer. Math. Soc.}, 353(2):557--584 (electronic), 2001.

\bibitem{Schultens2}
Jennifer Schultens.
\newblock The classification of {H}eegaard splittings for (compact orientable
  surface){$\,\times\, S\sp 1$}.
\newblock {\em Proc. London Math. Soc. (3)}, 67(2):425--448, 1993.

\bibitem{Schultens}
Jennifer Schultens.
\newblock Additivity of tunnel number for small knots.
\newblock {\em Comment. Math. Helv.}, 75(3):353--367, 2000.

\bibitem{Singer}
James Singer.
\newblock Three-dimensional manifolds and their {H}eegaard diagrams.
\newblock {\em Trans. Amer. Math. Soc.}, 35(1):88--111, 1933.

\bibitem{Souto}
Juan Souto.
\newblock Distances in the curve complex and the {H}eegaard genus.
\newblock Preprint. Available at {\tt
  www.picard.ups-tlse.fr/\%7Esouto/Heeg-genus.pdf}.

\end{thebibliography}

\end{document}